\documentclass[11pt,draft]{article} 
\usepackage{amsfonts,amsmath,amssymb,amsthm}

\newcommand{\rset}{\mathbf{R}}

\newcommand{\nset}{\mathbf{N}}
\newcommand{\ep}{\varepsilon}
\newcommand{\ind}{\mathbf{1}}
\newcommand{\fl}{\longrightarrow}
\newcommand{\e}{\mathbb{E}}
\newcommand{\p}{\mathbb{P}}
\newcommand{\lp}{\mathrm{L}}
\newcommand{\m}{\mathcal}

\newcommand{\ys}[1][p]{\mathcal{S}^{#1}}
\newcommand{\zs}[1][p]{\mathrm{M}^{#1}}
\usepackage{comment}
\usepackage[a4paper,body={160mm,230mm},centering]{geometry}
\newtheorem{thm}{Theorem}[section]
\newtheorem{lemme}[thm]{Lemma}
\newtheorem{prop}[thm]{Proposition}
\newtheorem{cor}[thm]{Corollary}
\theoremstyle{definition}
\newtheorem{df}[thm]{Definition}
\newtheorem{hyp}{Assumption}

\theoremstyle{remark}
\newtheorem{rem}[thm]{Remark}
\def\ds{\begin{displaystyle}}
\def\eds{\end{displaystyle}}
\def\dis{\displaystyle }
\def\<{\langle }
\def\>{\rangle }  
\newcommand{\nl}{\nolimits}

\def\R{\mathbf{ R}}
\def\N{\mathbf{N}}

\def\E{\mathbb E}

\def\P{\mathbb P}

\def\t{\tau}

\def\calb{{\cal B}}

\def\calf{{\cal F}}
\def\calg{{\cal G}}

\def\caln{{\cal N}}
\def\calp{{\cal P}}

\def\call{{\cal L}}

\title{\bf BSDEs with stochastic Lipschitz condition and quadratic PDEs in Hilbert spaces}

\author{Philippe Briand\\[.3em]
\normalsize\it
IRMAR, Universit\'{e} Rennes 1, 35042 Rennes Cedex, FRANCE\\
\normalsize\tt
philippe.briand@univ-rennes1.fr
\and
Fulvia Confortola\\[.3em]
\normalsize\it
Dipartimento di Matematica, Politecnico di Milano\\
\normalsize\it
piazza Leonardo da Vinci 32, 20133 Milano, Italy\\
\normalsize\tt
fulvia.confortola@unimib.it}

\date{January 29, 2006}

\begin{document}

\maketitle

\begin{abstract}
This paper is devoted to the study of the differentiability of solutions to real-valued backward stochastic differential equations (BSDEs for short) with quadratic generators driven by a cylindrical Wiener process. The main novelty of this problem consists in the fact that the gradient equation of a quadratic BSDE has generators which satisfy stochastic Lipschitz conditions involving BMO martingales. 
We show some applications to the nonlinear Kolmogorov equations.
\end{abstract}


\section{Introduction}

In this paper we are concerned with a real valued BSDE
\begin{equation*}
Y_\tau=  \Phi(X_T) + \int_\tau^T
  F (r ,X_r, Y_r,Z_r) \,dr -\int_t^T Z_r\, dW_r,
\qquad \tau\in [t,T],
\end{equation*}
where $W$ is a cylindrical Wiener process in some infinite dimensional Hilbert space $\Xi$ and the generator $F$ has quadratic growth with respect to the variable $z$. 
Quadratic BSDEs has been intensively studied by Kobylanski \cite{Kob00}, and then by Lepeltier and San Martin in \cite{LsM98} and more recently by Briand and Hu in \cite{BH06}. The process $X$, appearing in the generator and in the terminal value of the BSDE, takes its values in an an Hilbert space $H$ and it is solution of the following forward equation 
\begin{equation*}
  \left\{\begin{array}{l}\dis
dX_\tau=AX_\tau \, d\tau+
b(\tau,X_\tau)\, d\tau
+\sigma(\tau,X_\tau)\, dW_\tau,\quad \tau\in [t,T],
\\
\dis X_t=x\in H.
\end{array}
\right.
\end{equation*}
$A$ is the generator of a strongly continuous semigroup of
bounded linear operators $\{e^{tA}\}$ in $H$, $b$ and $\sigma$ are
functions with values in $H$ and $L_2(\Xi,H)$ -- the space of Hilbert-Schmidt operators from $\Xi$ to $H$ --  respectively.
Under suitable assumptions on the coefficients, there exists a unique adapted process
$(X,Y,Z)$ in the space $H\times \R \times L_2(\Xi,\R)$ solution to this forward-backward system.
The processes $X,Y,Z$ depend on the values of $x$ and $t$ occurring as initial conditions in the forward equation: we may denote them by $X^{t,x}$, $Y^{t,x}$ and $Z^{t,x}$.

Nonlinear BSDEs were first introduced by Pardoux and Peng \cite{PP90} and, since then, have been studied with great interest in finite and infinite dimensions: we refer the reader to \cite{EPQ97}, \cite{ElK97} and \cite{Par99} for an exposition of this subject and to \cite{MY99} for coupled forward-backward
systems. The interest in BSDEs comes from their connections with different mathematical fields, such as finance, stochastic control and partial differential equations. In this paper, we are concerned with the relation between BSDEs and nonlinear PDEs known as the nonlinear Feynman-Kac formula. More precisely, let us consider the following nonlinear PDE
$$
\partial_t u(t,x) +\call_t [u(t,\cdot)](x) +
F (t, x,u(t,x),\sigma(t,x)^*\nabla_xu(t,x)) = 0,\qquad u(T,x)=\Phi(x),
$$
where $ \call_t$ is the infinitesimal generator of the diffusion $X$. Then the solution $u$ is given by the formula $u(t,x)=Y^{t,x}_t$ which generalizes the Feynman-Kac formula to a nonlinear setting. 

Numerous results (for instance \cite{Pen91, PP92, Par98, Par99, Kob00}) show the connections between BSDEs set from a forward-backward system and solutions of a large class of quasilinear parabolic and elliptic PDEs. In the finite dimensional case, solutions to PDEs are usually understood in the viscosity sense. Here
we work in infinite dimensional spaces and consider solutions in the so called mild sense (see e.g. \cite{FT02}), which are intermediate between classical and viscosity solutions. This notion of solution seems natural in infinite dimensional framework: to have a mild solution its enough to prove that it is G\^{a}teaux differentiable. Hence we don't have to impose heavy assumptions on the coefficients as for the classical solutions. However a mild solution is G\^{a}teaux differentiable and thus more regular than a viscosity solution. For the probabilistic approach, this means that, in the infinite dimensional case, one has to study the regularity of  $X^{t,x}$, $Y^{t,x}$ and $Z^{t,x}$ with respect to $t$ and $x$ in order to solve the PDE.

This problem of regular dependence of the solution of a stochastic forward-backward system has been studied
in finite dimension by Pardoux, Peng \cite{PP92} and by El Karoui, Peng and Quenez \cite{EPQ97},
and, in infinite dimension, by Fuhrman and Tessitore  in \cite{FT02}, \cite{FT04}. In both cases, $F$ is assumed to be Lipschitz continuous with respect to $y$ and $z$.
In \cite{BC06P}, in infinite dimension, the generator $F$ is assumed to be only Lipschitz continuous only with respect to $z$ and monotone with respect to $y$ in the spirit of the works \cite{Pen91}, \cite{Par98} and more recently \cite{BDHPS03}.

In this work, we want to achieve this program when $F$ is quadratic with respect to $z$ meaning that the PDE is quadratic in the gradient. We will only consider the case of a bounded function $\Phi$. The study of the differentiability of the process $Y$ with respect to $x$  in this quadratic framework open an interesting problem of solvability of linear BSDEs with stochastic Lipschitz condition. Let us show with an example what happens in order to motivate the assumptions we will work with.

Let $(Y^x,Z^x)$ be the solution to the BSDE -- all processes are real in this example --
$$
Y_t^x = \Phi(x+W_t) + \frac{1}{2}\, \int_t^T \left| Z^x_s\right|^2 ds - \int_t^T Z^x_s\,dW_s
$$
where $\Phi$ is bounded and $\m C^1$. If $(G^x,H^x)$ stands for the gradient with respect to $x$ of $(Y^x,Z^x)$ then we have, at least formally,
$$
G^x_t = \Phi'(x+W_t) + \int_t^T Z^x_s H^x_s\, ds - \int_t^T H^x_s\, dW_s.
$$
In this linear equation, of course, the process $Z^x$ is not bounded in general so the usual Lipschitz assumption is not satisfied. It is only known that the process $Z^x$ is such that $\displaystyle \int_0^t Z^x_s\, dW_s$ is a BMO--martingale: this fact was used in \cite{HIM05} to prove a uniqueness result. BSDEs under stochastic Lipschitz condition have already been studied in \cite{EH97} and more recently in \cite{CRW06P}. However, the results in these papers do not fit our BMO-framework. This is the starting point of this paper.

The plan of the paper is as follows: Section \ref{sec-Notation} is devoted to notations. In Section \ref{secbsde} we recall some known results about BMO-martingales and we state a result of existence and uniqueness for BSDEs with generators satisfying a stochastic Lipschitz condition with  BMO feature. In section \ref{secfb} we apply the previous result to the study the regularity of the map $(t,x) \mapsto (Y_{\cdot}^{t,x},Z{\cdot}^{t,x})$ solution of the forward-backward system. The last section contain the applications to nonlinear Kolmogorov PDEs.

\section{Notations}
\label{sec-Notation}

\subsection{Vector spaces and stochastic processes}

In the following, all stochastic processes will be defined on
subsets of a fixed time interval $[0,T]$.

The letters $\Xi$, $H$ and $K$ will always denote Hilbert spaces.
Scalar product is denoted $\<\cdot,\cdot\>$, with a subscript
to specify the space if necessary. All
Hilbert spaces are assumed to be real and separable. $L_2(\Xi,K)$
is the space of Hilbert-Schmidt operators from $\Xi$ to $K$
endowed with the Hilbert-Schmidt norm. We observe that if $K=\R$ the space $L_2(\Xi, \rset)$ is the space $L(\Xi, \rset)$ of bounded linear operators from $\Xi$ to $\R$. By the Riesz isometry the dual space ${\Xi}^*=L(\Xi, \rset)$ can be identified with $\Xi$.

$W= \{W_t\}_{t\geq 0}$ is a cylindrical Wiener process with values in
  the infinite dimensional Hilbert space $\Xi$, defined on a probability
  space $(\Omega, \calf,\P)$; this means that a family
$W(t)$, $t\geq 0$, is a family of linear mappings from $\Xi$ to $L^2(\Omega)$ such that
\begin{description}
  \item[(i)] for every $u\in \Xi$, $\{ W(t)u,\; t\geq 0\}$ is a real
  (continuous) Wiener process;
  \item[(ii)] for every $u,v\in \Xi$ and $t\geq 0$,
  $\E\; (W(t)u\cdot W(t)v) = \<u,v\>_\Xi$.
\end{description}

  $\left\{\calf_t\right\}_{t\in [0,T]}$ will denote the natural
  filtration of $W$, augmented with the family $\caln$ of
  $\P$-null sets of $\calf_T$:
$$
\calf_t=\sigma(W(s)\; :\; s\in [0,t])\vee\caln.
$$
The filtration
$\left\{\calf_t\right\}_{t\in[0,T]}$ satisfies the usual conditions. All the concepts of
measurability for stochastic processes (e.g. predictability etc.)
refer to this filtration.
By $\calp$ we denote the
 predictable $\sigma$-algebra on $\Omega\times [0,T]$
and by $\calb(\Lambda)$ the Borel $\sigma$-algebra of any topological
space $\Lambda$.

Next we define several classes of stochastic processes which we use in the sequel.
For any real $p>0$, $\ys(K)$, or $\ys$ when no confusion is possible, denotes the set of $K$-valued, adapted and c\`adl\`ag processes
$\{Y_t\}_{t\in[0,T]}$ such that
$$
\left\| Y \right\|_{\ys} :=\e\left[\sup\nl_{t\in[0,T]} |Y_t|^p\right]^{1\wedge 1/p} < +\infty.
$$
If $p\geq 1$, $\|\cdot\|_{\mathcal{S}^p}$ is a norm on
$\mathcal{S}^p$ and if
$p\in(0,1)$, $(X,X')\longmapsto \big\|X-X'\big\|_{\mathcal{S}^p}$ defines a distance on $\ys$. Under this metric, $\ys$ is
complete. $\zs$ ($\zs\left(L_2(\Xi, K)\right)$) denotes the set of (equivalent classes of)
predictable processes $\{Z_t\}_{t\in[0,T]}$ with values
in $L_2(\Xi, K)$ such that
$$
\left\| Z \right\|_{\zs} : = \e\left[\Big(\int_0^T |Z_s|^2\,
ds\Big)^{p/2}\right]^{1 \wedge 1/p} < +\infty.
$$
For $p\geq 1$, $\zs$ is a Banach space endowed with
this norm  and for $p\in(0,1)$, $\zs$ is a complete metric
space with the resulting distance. We set $\ys[] = \cup_{p>1} \ys$, $\zs[] = \cup_{p>1} \zs$ and $\ys[\infty]$ stands for the set of predictable bounded processes.

Given an element $\Psi$ of
$L_{\calp}^2(\Omega\times [0,T];L_2(\Xi,K))$, one can define the It\^o stochastic integral
$\int_0^t\Psi(\sigma)\, dW_{\sigma}$, $t\in [0,T]$;
it is a $K$-valued martingale with continuous path such that
$$ \e\left[\sup\nl_{t\in[0,T]} |\int_0^t\Psi(\sigma)\, dW_{\sigma}|^2\right]^{ 1/2} < +\infty.$$

The previous definitions have obvious extensions to processes
defined on subintervals of $[0,T]$.

\subsection{The class $\calg $}

 $F:X\to V$, where $X$ and $V$ are two Banach spaces, has  a directional derivative at point $x\in X$ in the direction $h\in X$ when
  $$
  \nabla F(x;h)=\lim_{s\to 0}\frac{F(x+sh)-F(x)}{s},
  $$
 exists in the topology of $V$. $F$ is said to be G\^ateaux
differentiable at point $x$ if $\nabla F(x;h)$  exists for every $h$ and
there exists an element of $L(X,V)$, denoted $\nabla F(x)$ and called G\^ateaux
derivative, such that $\nabla F(x;h)=\nabla F(x)h$ for every $h\in X$.

\begin{df} 
     $F:X\to V$ belongs to the class $\calg^1 (X;V)$ if it is continuous,
G\^ateaux differentiable on $X$, and  $\nabla F:X\to L(X,V)$
is strongly continuous.
\end{df}
 
In particular, for every $h\in X$ the map $\nabla F(\cdot)h:X\to V$ is continuous. Let us recall some features of the class $\calg^1(X,V)$ proved in \cite{FT02}.

\begin{lemme}\label{proprietadig}
  Suppose $F\in\calg^1(X,V)$. Then
\begin{description}
  \item[(i)] $(x, h)\mapsto \nabla F(x)h$ is
  continuous from $X\times X$ to $V$;
  \item[(ii)]  If $G\in
  \calg^1(V,Z)$ then $G( F)\in \calg^{1} (X,Z)$ and
  $\nabla(G( F))(x)=\nabla G(F(x))\nabla F(x)$.
\end{description}
\end{lemme}

\begin{lemme}
\label{staing}
  A map $F:X\to V$ belongs to $\calg^1(X,V)$
  provided the following conditions hold:
\begin{description}
  \item[(i)]  the directional
  derivatives $\nabla F(x;h)$ exist
  at every point $x\in X$ and in every direction $h\in X$;
  \item[(ii)] for every  $h$, the mapping
   $\nabla F(\cdot;h): X\to V$
  is continuous;
  \item[(iii)] for every $x$, the mapping $h\mapsto\nabla F(x;h)$
  is continuous from $X$ to $V$.
\end{description}
\end{lemme}

These definitions can be generalized to functions depending on
several variables.  For instance, if $F$ is a function from $X\times Y$ into $V$, the partial
directional and G\^ateaux derivatives with respect to the first argument, at point $(x,y)$ and in the
direction  $h\in X$, are denoted $\nabla_x F(x,y;h)$ and $\nabla_x F(x,y)$ respectively.

\begin{df}
$F:X\times Y\to V$ belongs to the class $\calg^{1,0} (X\times Y;V)$ if it is continuous,
G\^ateaux differentiable with respect to
$x$ on $X\times Y$, and  $\nabla_xF:X\times Y\to L(X,V)$
is strongly continuous.
\end{df}

As in Lemma \ref{proprietadig}, the map
$(x, y, h)\mapsto \nabla_xF(x,y)h$ is
  continuous from $X\times Y\times X$ to $V$, and
the  chain rules hold. One can also extend Lemma
\ref{staing} in the following way.

\begin{lemme}\label{staing2}
  A continuous map $F:X\times Y\to V$ belongs to $\calg^{1,0}(X\times Y,V)$ provided the following conditions hold:
\begin{description}
  \item[(i)]
  the directional
  derivatives $\nabla_x F(x,y;h)$ exist
  at every point $(x,y)\in X\times Y$
  and in every direction $h\in X$;
  \item[(ii)] for every  $h$, the mapping
   $\nabla F(\cdot,\cdot ;h): X\times Y\to V$
  is continuous;
  \item[(iii)] for every $(x,y)$,
   the mapping $h\mapsto\nabla_x F(x,y;h)$
  is continuous from $X$ to $V$.
\end{description}
\end{lemme}

When $F$ depends on additional arguments,
 the previous definitions and properties  have
obvious generalizations. For instance, we say that
$F:X\times Y\times Z\to V$
belongs to  $\calg^{1,1,0} (X\times Y\times Z;V)$
if it is continuous,
G\^ateaux differentiable with respect to
$x$ and $y$ on $X\times Y\times Z$,
and  $\nabla_xF:X\times Y\times Z\to L(X,V)$ and
$\nabla_yF:X\times Y\times Z\to L(Y,V)$ are strongly continuous.

\section{BSDEs with random Lipschitz condition}
\label{secbsde}
In this section, we want to study the  BSDE
\begin{equation}
\label{eqgen}
Y_t = \xi + \int_t^T f(s,Y_s,Z_s)\, ds - \int_t^T Z_s\, dW_s
\end{equation}
when the generator $f$ is Lipschitz but with random Lipschitz constants. This kind of BSDEs were also considered in \cite{EH97} and more recently in \cite{CRW06P}. However our framework is different from the setting of the results obtained in these papers. Let us recall that a generator is a random function
$f:[0,T]\times\Omega \times \R \times\ L_2(\Xi, \R)
\fl \R$ which is measurable with respect to $\m P\otimes\mathcal{B}(\R)\otimes\mathcal{B}(\Xi)$
and a terminal condition is simply a real $\mathcal{F}_T$--measurable random variable.
From now on, we deal only with generators such that, $\p$--a.s., for each $t\in[0,T]$, $(y,z)\fl f(t,y,z)$
is continuous.

By a solution to the BSDE~\eqref{eqgen} we mean a pair $(Y,Z)=\{ (Y_t,Z_t)\}_{t\in[0,T]}$ of predictable processes  with values in $\rset\times L_2(\Xi, \rset)$
such that $\p$--a.s., $t\longmapsto Y_t$ is continuous, $t\longmapsto Z_t$ belongs to  $\lp^2(0,T)$, $t\longmapsto f(t,Y_t,Z_t)$ belongs to $\lp^1(0,T)$
and $\p$--a.s.
$$
Y_t = \xi +\int_t^T f(s,Y_s,Z_s)\, ds -\int_t^T Z_s\, dW_s, \qquad
0\leq t\leq T.
$$

We will work with the following assumption on the generator.  

\begin{hyp}
\label{hmon}
There exist a real process $K$ and a constant $\alpha\in(0,1)$  such that $\p$--a.s.:
\begin{itemize} 
\item for each $t\in[0,T]$, $(y,z)\fl f(t,y,z)$ is continuous ;
\item for each $(t,z)\in[0,T]\times L_2(\Xi, \R)$,
$$ \forall  y,p \in \R,\qquad (y-p)(f(t,y,z) -f(t,p,z)) \leq K_t^{2\alpha}\, |y-p|^2$$

\item for each $(t,y)\in[0,T]\times \R$,
$$
\forall \left(z,q\right)\in L_2(\Xi,\R)\times L_2(\Xi, \R),\qquad  \left| f(t,y,z) -f(t,y,q)\right| \leq K_t\, |z-q|_{L_2(\Xi,\R)}.
$$
\end{itemize}
\end{hyp} 

In the classical theory, the process $K$ is constant but for the application we have in mind we will only assume the following.  

\begin{hyp}
\label{hbmo}
$\{ K_s \}_{s\in[0,T]}$ is a predictable real process  bounded from below by 1 such that there is a constant $C$ such that, for any stopping time $\tau\leq T$,
$$
\e\left( \int_\tau^T |K_s|^2\,ds \:\Big | \: \m F_\tau  \right) \leq C^2.
$$
$N$ denotes the smallest constant $C$ for which the previous statement is true.
\end{hyp}  

This assumption says that, for any $u\in L_2(\Xi,\R)$ such that $||u||_{L_2(\Xi,\R)}=1$ the martingale 
$$
M_t = \int_0^t K_s u dW_s,\quad 0\leq t\leq T
$$
is a BMO-martingale with $\| M \|_{BMO_2} = N$. We refer to \cite{Kaz94} for the theory of BMO--martingales and we just recall the properties we will use in the sequel. It follows from the inequality (\cite[p. 26]{Kaz94}),
$$
\forall n\in\nset^*,\qquad \e\left[\langle M \rangle_T ^n \right]= \e\left[\Big( \int_0^T |K_s|^2\, ds \Big)^n \right] \leq n!\, N^{2n}
$$
that $M$ belongs to $\mathrm{H}^p$ for all $p\geq 1$ and moreover
\begin{equation}
\label{eqexpmom}
\forall \alpha\in(0,1),\quad \forall p\geq 1,\qquad \eta(p)^p:=\e\left[\exp\left(p\int_0^T |K_s|^{2\alpha} ds\right)\right] < +\infty.
\end{equation} 
The very important feature of BMO--martingales is the following: the exponential martingale
$$
\m E(M)_t=\m E_t = \exp\left(\int_{0}^t K_{s} u\cdot dW_{s} - \frac{1}{2} \int_{0}^t |K_{s}|^2 ds\right)
$$
is a uniformly integrable martingale. More precisely, $\{ \m E_t \}_{0\leq t\leq T}$ satisfies a reverse H\"older inequality. Let $\Phi$ be the function defined on $(1,+\infty)$ by
$$
\Phi(p) = \left( 1 + \frac{1}{p^2} \log \frac{2p-1}{2(p-1)}\right)^{1/2}-1~;
$$
$\Phi$ is nonincreasing with $\lim_{p\to 1} \Phi(p) = +\infty$, $\lim_{p\to +\infty} \Phi(p) = 0$. Let $q_*$ be such that $\Phi(q_*) = N$. Then, for each $1<q<q_*$ and  for all stopping time $\tau\leq T$,
\begin{equation}
\label{eqrevhol}
\e\left( \m E(M)_T^{q} \: \big|\: \m F_\tau \right) \leq K(q,N)\, \m E(M)_\tau^{q}
\end{equation}
where the constant $K(q,N)$ can be chosen depending only on $q$ and $N=\| M \|_{BMO_2}$ e.g.
$$
K(q,N) = \frac{2}{1-2(q-1)(2q-1)^{-1}\exp(q^2(N^2+2N))}.
$$

\begin{rem}
\label{restriv}  
If we denote $\p^*$ the probability measure on $(\Omega, \m F_T)$ whose density with respect to $\p$ is given by $\m E_T$ then $\p$ and $\p^*$ are equivalent.

Moreover, it follows from \eqref{eqrevhol} and H\"older's inequality that, if $X$ belongs to $\lp^p(\p)$ then $X$ belongs to $\lp^s(\p^*)$ for all $s<p/p_*$ where $p_*$ is the conjugate exponent of $q^*$.
\end{rem}

We assume also some integrability conditions on the data.
For this, let $p_*$ be the conjugate exponent of $q_*$.

\begin{hyp}
\label{hint}
There exists $p^* > p_*$ such that
$$
 \e\left[ |\xi|^{p^*} + \Big(\int_{0}^T |f(s,0,0)|\,ds \Big)^{p^*}\right] < +\infty.
$$
\end{hyp} 

As usual for BSDEs, we begin with some apriori estimate. The first one shows that, one can control the process $Y$ as soon as the process $Z$ has some integrability property. The following lemma relies heavily on the reverse H\"older's inequality.

\begin{lemme}
\label{resz2y}
Let the assumptions \ref{hmon}, \ref{hbmo} and \ref{hint} hold.
If $(Y,Z)$ is a solution to \eqref{eqgen} such that, for some $r>p_*$, $Z\in\zs[r]$,  then, for each $p\in(p_*,p^*)$, $Y\in\ys$ and
$$
\left\| Y \right\|_{\ys} \leq C \left\| |\xi| + \int_{0}^T |f(s,0,0)|\,ds  \right\|_{p^*},
$$
for a suitable constant $C$ depending on $p$, $p^*$, $p_*$ and $N$.
\end{lemme}

\begin{proof}
The starting point to obtain this estimate is a linearization of the generator of the BSDE~\eqref{eqgen}. Let us set
$$
a_s = \frac{f(s,Y_s,Z_s)-f(s,0,Z_s)}{Y_s},\qquad b_s=\frac{f(s,0,Z_s)-f(s,0,0)}{|Z_s|_{L_2(\Xi,\R)}^2} Z_s.
$$
Then, $(Y,Z)$ solves the linear BSDE
$$ Y_t = \xi + \int_t^T \left( f(s,0,0) + a_s\, Y_s + <b_s, Z_s>_{L_2(\Xi,\R)} \right)  ds - \int_t^T Z_s\, dW_s.$$

As usual, let us set $e_t = e^{\int_0^t a_s\, ds}$. We have,
$$
e_t Y_t = e_T \xi + \int_t^T e_s f(s,0,0)\,ds - \int_t^T e_sZ_s\cdot dW^*_s,
$$
where we have set $W^*_s = W_s - \int_0^s b_r\, dr$.  Of course, we want to take the conditional expectation of the previous equality with respect to the probability $\p^*$ whose density is 
$$
\m E(I(b))_T=\exp\left(\int_0^T b_s dW_s -\frac{1}{2} \int_0^T |b_s|_{L_2(\Xi,\R)}^2 ds\right)
$$ under which $B^*$ is a Brownian motion. To do this, let us observe that $|b_s|_{L_2(\Xi,\R)} \leq K_s$ so that $\| I(b) \|_{BMO_2} \leq \| M \|_{BMO_2}$ and $\m E(I(b))$ satisfies the reverse H\"older inequality \eqref{eqrevhol} for all $q<q_*$ (with the same constant).

Moreover, it follows from \ref{hmon} that $ a_s \leq K_s^{2\alpha}$ and, in particular, \eqref{eqexpmom} says that the process $e$ belongs to all $\ys$ spaces. Thus $e_T\xi$ belongs to $\lp^{p}$ for all $p<p_*$ and the same is true for $\int_0^T e_s |f(s,0,0)|\,ds$. In the same way, we have, for all $\rho<r$,
$$
\e\left[ \Big( \int_0^T e_s^2 |Z_s|^2 ds \Big)^{\rho/2} \right] \leq \e\left[ \sup e_t^\rho \Big( \int_0^T |Z_s|^2 ds \Big)^{\rho/2}\right] < +\infty.
$$

Using Lemma~\ref{restriv}, we deduce that $e_T\xi$ and $\int_0^T e_s |f(s,0,0)|\,ds$ belongs to $\lp^p(\p^*)$ for all $p<p^*/p_*$ and $\Big(\int_0^T |Z_s|^2 ds \Big)^{1/2}$ belongs to $\lp^s$ for all $s<r/p_*$.

Thus we can take the conditional expectation to obtain
$$
e_tY_t = \e^*\left( e_T \xi + \int_t^T e_s f(s,0,0)\,ds \: \Big|\: \m F_t\right),
$$
and, as a byproduct of this equality,  we get
$$
|Y_t| \leq (\m E_t)^{-1} \e\left( \m E_T \left( |\xi| e_T/e_t  + \int_t^T |f(s,0,0)| e_s/e_t\,ds \right) \: \Big|\: \m F_t\right).
$$
Taking into account \ref{hmon}, we have $a_s \leq K_s^{2\alpha}$ and, for all $s>t$,
$$
e_s/e_t \leq \exp\left(\int_t^s K_r^{2\alpha}\,dr\right) \leq \exp\left(\int_0^T K_r^{2\alpha}\, dr\right),
$$
from which we deduce the inequality
$$
|Y_t| \leq (\m E_t)^{-1} \e\left( \m E_T \Gamma_T X \: \big|\: \m F_t\right),
$$
where we have set
$$
\Gamma_T=\exp\left(\int_0^T K_r^{2\alpha}\, dr\right),\quad\text{and}\quad X= \left( |\xi| + \int_0^T |f(s,0,0)|\,ds \right).
$$
Using the reverse H\"older inequality, for each $r>p_*$, we have, $q=r/(r-1) < q_*$ and
$$
|Y_t| \leq (\m E_t)^{-1} \e\left( \m E_T ^q  \: \big|\: \m F_t\right)^{1/q} \e\left( \Gamma_T^r X^r \: \big|\: \m F_t\right)^{1/r} \leq K(q,N)^{1/q} \e\left( \Gamma_T^r X^r \: \big|\: \m F_t\right)^{1/r}
$$
Doob's inequality gives for all $p_*<r<p$,
$$
\e\left[ \sup_{t\in[0,T]} |Y_t|^p \right] \leq K(q,N)^{p/q} \left(\frac{p}{p-r}\right)^{p/r} \e[\Gamma_T^p X^p].
$$
Now, let $p\in(p_*,p^*)$, from H\"older inequality, we have, for each $p_*<r<p$,
$$
\e\left[ \sup\nl_{t\in[0,T]} |Y_t|^p \right] \leq K(q,N)^{p/q} \left(\frac{p}{p-r}\right)^{p/r} \eta\left(pp^*/(p^*-p)\right)^{p}\e[X^{p^*}]^{p/p^*} .
$$
It follows that, for $p_*<r<p<p^*$,
$$
\| Y \| _{\ys} \leq K\left( \frac{r}{r-1}, N\right)^{(r-1)/r} \left(\frac{p}{p-r}\right)^{1/r} \eta\left(\frac{pp^*}{p^*-p} \right) \left\| |\xi| + \int_0^T |f(s,0,0)|\,ds \right\|_{p^*},
$$
which gives the result taking $r=(p+p_*)/2$.
\end{proof}

We keep on by showing that on can obtain an estimate for the process $Z$ in terms of the norm of $Y$. This kind of results is quite classical see e.g. \cite{BDHPS03}. We give the proof in our framework for  the ease of the reader.

\begin{lemme}
\label{resy2z}
Let us assume that
$$
y\cdot f(t,y,z) \leq |y| f_t + K_t^{2\alpha} |y|^2 + K_t |y|\,|z|
$$
for nonnegative processes $f$ and $K$.

If $(Y,Z)$ solves the BSDE~\eqref{eqgen}, with $Y\in\ys[q]$ then, for each $p<q$, $Z\in\zs$ and
$$
\left\| Z \right\|_{\zs} \leq C\left( \| Y \|_{\ys} + \left\|\int_0^T f_s\, ds \right\|_p + \| Y \|_{\ys[q]} \left\|\Big(\int_0^T \left( K_s^{2\alpha} + K_s^2 \right) ds\Big)^{1/2} \right\|_{pq/(q-p)} \right),
$$
where $C$ depends only on $p$ and $q$.
\end{lemme}

\begin{proof}
We follow \cite{BDHPS03}. For each integer $n\geq 1$, let us introduce the stopping time
$$
\tau_n =\inf\left\{ t\in[0,T],\: \int_0^t |Z_r|^2\, dr \geq
n\right\}\wedge T.
$$
It\^o's formula gives us,
$$
|Y_0|^2 + \int_0^{\tau_n} |Z_r|^2\, dr = |Y_{\tau_n}|^2 +
2 \int_0^{\tau_n} \langle Y_r, f(r,Y_r,Z_r)\rangle\,
dr-2\int_0^{\tau_n}\langle Y_r,Z_rdW_r\rangle.
$$
But, from the assumption on $f$, we have,
$$
2y\cdot f(r,y,z) \leq 2|y|f_r +2K_r^{2\alpha}|y|^2 + 2K_r^2|y|^2 +
|z|^2/2.
$$
Thus, since $\tau_n\leq T$, we deduce that
$$
\frac{1}{2}\int_0^{\tau_n} |Z_r|^2\, dr  \leq   Y_*^2 + 2 Y_*
\int_0^T f_r\, dr + 2 Y_*^2 \int_0^T \left( K_r^{2\alpha}  + K_r^2 \right) dr+ 2\Big|\int_0^{\tau_n} \langle Y_r,Z_r dW_r\rangle
\Big|.
$$
It follows that
$$
\int_0^{\tau_n} |Z_r|^2\, dr \leq 4\left( Y_*^2 + \Big(\int_0^T
f_r\, dr\Big)^2 + Y_*^2 \int_0^T \left( K_r^{2\alpha}  + K_r^2 \right) dr+ \Big|\int_0^{\tau_n} \langle Y_r,Z_r
dW_r\rangle\Big|\right)
$$
and thus that
\begin{equation}
\label{conclure}
\begin{split}
\lefteqn{\Big(\int_0^{\tau_n} |Z_r|^2\, dr\Big)^{p/2} } \qquad   \\
  & \leq c_p \left(Y_*^p + \Big(\int_0^T
f_r\, dr\Big)^p + Y_*^{p} \Big( \int_0^T \left( K_r^{2\alpha}  + K_r^2 \right) dr \Big)^{p/2}+\Big|\int_0^{\tau_n}
\langle Y_r,Z_r dW_r\rangle\Big|^{p/2}\right).
\end{split}
\end{equation}

But by the BDG inequality, we get
$$
c_p\,\e\left[\Big|\int_0^{\tau_n}\! \langle Y_r,Z_r
dW_r\rangle\Big|^{p/2}\right]\leq d_p\,\e\left[\left(\int_0^{\tau_n}
\! |Y_r|^2\,|Z_r|^2 dr\right)^{p/4}\right] \leq
d_p\,\e\left[Y_*^{p/2}\Big(\int_0^{\tau_n}\!
|Z_r|^2 dr\Big)^{p/4}\right],
$$
and thus
$$
c_p\,\e\left[\Big|\int_0^{\tau_n} \langle Y_r,Z_r dW_r\rangle\Big|^{p/2}\right]
\leq  \frac{d_p^2}{2} \e\left[Y_*^p\right]
+\frac{1}{2}\, \e\left[\Big(\int_0^{\tau_n} |Z_r|^2\,
dr\Big)^{p/2}\right].
$$

Coming back to the estimate~(\ref{conclure}), we get, for each $n\geq 1$,
$$
\e\left[\Big(\int_0^{\tau_n} |Z_r|^2\,
dr\Big)^{p/2}\right] \leq C_p\,\e\left[Y_*^p + \Big(\int_0^T f_r\,
dr\Big)^p + Y_*^{p} \Big( \int_0^T \left( K_s^{2\alpha}  + K_s^2 \right) ds \Big)^{p/2}\right]
$$
and, Fatou's lemma implies that
$$
\e\left[\Big(\int_0^T |Z_r|^2\,
dr\Big)^{p/2}\right] \leq C_p\,\e\left[ Y_*^p + \Big(\int_0^T f_r\,
dr\Big)^p + Y_*^{p} \Big( \int_0^T \left( K_s^{2\alpha} + K_s^2 \right) ds \Big)^{p/2}\right].
$$
The result follows from H\"older's inequality.
\end{proof}  

The previous two lemmas lead the following result.

\begin{cor}
\label{resutil}
Let the assumptions \ref{hmon}, \ref{hbmo} and \ref{hint} hold. If $(Y,Z)$ is a solution to \eqref{eqgen} such that, for some $r>p_*$, $Y\in\ys[r]$,  then, for each $p\in(p_*,p^*)$, $(Y,Z)\in\ys\times\zs$ and
$$
\left\| Y \right\|_{\ys} + \left\| Z \right\|_{\zs} \leq C\, \left\| |\xi| + \int_{0}^T |f(s,0,0)|\,ds  \right\|_{p^*} \left( 1 + \left\| \Big(\int_0^T \left(K_s^{2\alpha} + K_s^2 \right) ds\Big)^{1/2} \right\|_{p(p^*+p)/(p^*-p)} \right)
$$
where $C$ depends on $p$, $p_*$, $p^*$ and $N$.
\end{cor}

\begin{proof}
Since $Y$ belongs to $\ys$ for some $p>p_*$, there exists by Lemma~\ref{resy2z} $r\in(p_*,p^*)$ such that $Z$ belongs to $\zs[r]$. It follows from Lemma~\ref{resz2y} that $Y$ belongs to $\ys$ for all $p<p^*$ and then by Lemma~\ref{resy2z} $Z\in\zs$ for all $p<p^*$.

The inequality comes from the choice $q=(p+p^*)/2$ in Lemma~\ref{resy2z} together with the estimate of Lemma~\ref{resz2y}.
\end{proof}

\begin{hyp}
\label{hgro}
There exists a nonnegative predictable process $f$ such that, 
\begin{equation*}
    \e\left[\Big(\int_0^T f(s)\,ds\Big)^{p^*}\right] < +\infty
\end{equation*}
and $\p$--a.s.
$$
\forall (t,y,z)\in[0,T]\times\rset\times L_2(\Xi,\rset),\qquad
\left| f(t,y,z) \right| \leq f(t) + K_t^{2\alpha} |y| + K_t |z|.
$$
\end{hyp}

\begin{thm}
\label{resesun}
Let the assumptions \ref{hmon}, \ref{hbmo}, \ref{hint} and \ref{hgro} hold. Then BSDE~\eqref{eqgen} has a unique solution $(Y,Z)$ which belongs to $\ys\times\zs$ for all $p<p^*$.
\end{thm}

\begin{proof}
Let us prove first uniqueness.  Let $(Y^1,Z^1)$ and $(Y^2,Z^2)$ be solutions to \eqref{eqgen} such that $Y^1$ and $Y^2$ belongs to $\ys$ for $p>p_*$. The by Corollary~\ref{resutil}, $(Y^1,Z^1)$ and $(Y^2,Z^2)$ belongs to $\ys\times\zs$ for all $p<p_*$. Moreover, $U=Y^1-Y^2$ and $V=Z^1-Z^2$ solves the BSDE
    $$
    U_t = \int_t^T F(s,U_s,V_s)\,ds - \int_t^T V_s\cdot dW_s,
    $$
    where $F(t,u,v) = f\left(t,Y^2_t+u,Z^2_t+v\right) - f\left(t,Y^2_t,Z^2_t\right)$. We have $F(t,0,0)=0$ and $F$ satisfies \ref{hmon} with the same process $K$. It follows from Corollary~\ref{resutil} that $(U,V)\equiv (0,0)$.
    
    Let us turn to existence. For each integer $n\geq 1$, let $\tau_n$ be the following stopping time:
$$
\tau_n=\inf\left\{t\in[0,T] : \int_0^t \left( f(s) + K_s^2 \right)  ds  \geq n \right\} \wedge T.
$$
Let $\xi^n = \xi\ind_{|\xi|\leq n}$ and $(Y^n,Z^n)$ be the solution to the BSDE
$$
Y^n_t = \xi^n + \int_t^T \ind_{s\leq \tau_n} f\left(s,Y^n_s,Z^n_s\right)\, ds - \int_t^T Z^n_s  dW_s.
$$
The existence of the solution $(Y^n,Z^n)$ to the previous equation comes from \cite{Mor06P}. Indeed, we have, setting $f^n(t,y,z) =  \ind_{t\leq \tau_n} f(t,y,z)$,
$$
\left| f^n(t,y,z) \right| \leq \ind_{t\leq \tau_n} \left(f(t) + K_t^{2\alpha} + K_t^2/2 \right) \left(1+|y|\right) + |z|^2/2,
$$
and, $\p$--a.s.
$$
\int_0^T \ind_{t\leq \tau_n} \left(f(t) + K_t^{2\alpha} + K_t^2/2 \right) dt \leq 5n/2.
$$
Since $\xi^n$ is bounded by $n$, the previous BSDE has a unique solution $(Y^n,Z^n)$ such that $Y^n$ is a bounded process and $Z^n\in\zs[2]$. Since
$$
\int_0^T \left| f^n(t,0,0)\right|\,dt \leq n,
$$
we know, from Corollary~\ref{resutil}, that $(Y^n,Z^n)\in \ys\times\zs$ for all $p$.

Moreover, still by Corollary~\ref{resutil}, the sequence $\left( (Y^n,Z^n) \right)_{n\geq 1}$ is bounded in $\m K^p:=\ys\times\zs$ for all $p<p^*$.

Let us show that $\left( (Y^n,Z^n) \right)_{n\geq 1}$ is a Cauchy sequence in $\m K^p:=\ys\times\zs$ for all $p<p^*$. Let $m>n\geq 1$ and let us set as before $U = Y^m-Y^n$, $V=Z^m-Z^n$. Then $(U,V)$ solves the BSDE
$$
U_t = \xi^m-\xi^n + \int_t^T F(s,U_s,V_s) \, ds - \int_t^T V_s dW_s
$$
where
$$
F(t,u,v) = \ind_{t\leq \tau_m}\left( f\left(t,u+Y^n_t,v+Z^n_t\right)-f\left(t,Y^n_t,Z^n_t\right) \right) - \ind_{\tau_n<t\leq \tau_m} f\left(t,Y^n_t,Z^n_t\right).
$$
$F$ satisfies \ref{hmon} and $F(t,0,0)= - \ind_{\tau_n<t\leq \tau_m} f\left(t,Y^n_t,Z^n_t\right)$ belongs to $\lp^p$ for all $p\geq 1$.

Since $\xi\in\lp^{p^*}$, $\left\| \xi^m - \xi^n \right\|_{p^*} \fl 0$ if $n\to\infty$. Moreover, we have, from \ref{hgro} and H\"older inequality,
$$
\int_0^T \left| F(t,0,0) \right| dt \leq  \int_{\tau_n}^T f(t)\, dt + \sup\nl_t \left| Y^n_t\right| \, \int_{\tau_n}^T K_t^{2\alpha} \, dt + \left(\int_{\tau_n}^T K_t^2\,dt \right)^{1/2} \left(\int_0^T \left| Z^n_t \right|^2\,dt \right)^{1/2}.
$$

Let $p<p^*$. We choose $p<q<r<p^*$. It follows from the previous inequality, using H\"older inequality, that
$$
\left\| \int_0^T \left| F(t,0,0) \right| dt \right\|_{q} \leq \left\| \int_{\tau_n}^T f(t)\, dt\right\|_q + \left\| Y^n \right\|_{\ys[r]} \left\| \int_{\tau_n}^T K_t^{2\alpha} \, dt \right\|_{\frac{qr}{r-q}} + \left\| Z^n \right\|_{\zs[r]} \left\| \Big(\int_{\tau_n}^T K_t^{2} \, dt\Big)^{\frac{1}{2}} \right\|_{\frac{qr}{r-q}}.
$$
Let us recall that $\tau_n\to T$ $\p$--a.s and that the sequence $\left( (Y^n,Z^n) \right)_{n\geq 1}$ is bounded in $\m K^r$. Since $\int_0^T f(t)\, dt$ belongs to $\lp^{p^*}$, $\int_0^T K_t^{2\alpha}\, dt$ and $\int_0^T K_t^2\,dt$ has moments of all order, the right hand side of the previous inequality tends to 0 as $n$ tends to infinity.

It follows from Corollary~\ref{resutil} -- applied with $q$ instead of $p^*$ -- that $\left( (Y^n,Z^n) \right)_{n\geq 1}$  is a Cauchy sequence in $\m K^p$ and this is valid as soon as $p<p^*$.

It is easy to check that the limit of this sequence is a solution to BSDE~\eqref{eqgen}
\end{proof}

\section{The forward-backward system}
\label{secfb}

In this section, we apply the previous results on BSDEs to study the differentiability of the solution to the  following quadratic BSDE  
\begin{equation}\label{eq:markb}
 Y_{\tau}^{t,x} = \Phi\left(X_T^{t,x}\right) + \int_{\tau}^T F\left(r,X_r^{t,x},Y_r^{t,x},Z_r^{t,x}\right) dr -\int_{\tau}^T Z^{t,x}_r\,dWr, \quad 0\leq \tau\leq T,
\end{equation}
where $\left\{ X^{t,x}_\tau\right\}_{0\leq t\leq \tau}$  is the solution to 
\begin{equation}\label{eq:markf}
X^{t,x}_\tau = e^{(\tau-t)A}x+\int_t^\tau e^{(\tau-r)A}
b\left(r,X^{t,x}_r\right)\, dr
+\int_t^\tau e^{(\tau-r)A}
\sigma\left(r,X^{t,x}_r\right)\, dW_r,\quad t\leq \tau \leq T.
\end{equation} 
As usual, we have set $X^{t,x}_\tau=x$ for $\tau<t$.  Of course, from It\^o's formula, we have
$$
   dX^{t,x}_\tau = AX^{t,x}_\tau\, d\tau+b\left(\tau,X^{t,x}_\tau\right)\, d\tau
+\sigma\left(\tau,X^{t,x}_\tau\right)\, dW_\tau,\quad \tau\in [t,T],\qquad X^{t,x}_\tau=x\in H,\quad \tau\leq t.
$$ 
But a solution of this equation is always understood as an $(\calf_t)$-predictable continuous process $X$ solving~\eqref{eq:markf}.

We will work under the following assumption on the diffusion coefficients.

\begin{hyp}
\label{ipotesiuno}
\begin{description}
  \item[(i)] The operator $A$ is the generator of a strongly
  continuous semigroup $e^{tA}$, $t\geq 0$, in the Hilbert space
  $H$.
  \item[(ii)] The mapping $b:[0,T]\times H\to H$ is measurable
  and satisfies, for
  some constant $L>0$,
 \begin{eqnarray*}
  |b(t,x)-b(t,y)| & \leq & L\, |x-y|,\qquad t\in [0,T],\; x,y\in H, \\
|b(t,x)| & \leq & L\,   (1+|x|), \qquad t\in [0,T],\; x \in H.
\end{eqnarray*}
  \item[(iii)] $\sigma: [0,T]\times H \fl L(\Xi,H)$ is such that,
  for every $v\in \Xi$, the map $\sigma v: [0,T]\times H\to H$ is measurable,
  $e^{sA}\sigma(t,x)\in L_2(\Xi,H)$ for every $s>0$,
  $t\in [0,T]$ and $x\in H$,
   and
\begin{eqnarray*}
|e^{sA}\sigma(t,x)|_{L_2(\Xi,H)} & \leq & L\; s^{-\gamma}
  (1+|x|),\\
  |e^{sA}\sigma(t,x)-e^{sA}\sigma(t,y)|_{L_2(\Xi,H)} & \leq  & L\; s^{-\gamma}
  |x-y|, \\
|\sigma(t,x)|_{L(\Xi,H)} & \leq & L\;   (1+|x|),
\end{eqnarray*}
 for
  some constants $L>0$ and $\gamma\in [0,1/2)$.
\item[(iv)] For every $s>0$,
$t\in [0,T]$,
$$b(t,\cdot)
\in \calg^1(H, H),\qquad
e^{sA}\sigma(t,\cdot)\in\calg^1(H, L_2(\Xi,H)).
$$
  \end{description}
\end{hyp}

A consequence of the previous assumptions is that, for every $s>0$, $t\in [0,T]$, $x,h\in H$,
 $$
  |\nabla_x b(t,x)h|\leq L\; |h|, \qquad
|\nabla_x(e^{sA}\sigma(t,x))h|_{L_2(\Xi,H)}\leq L\; s^{-\gamma} |h|.
$$

The following results are proved by Fuhrman and Tessitore in \cite{FT02}.

\begin{prop}\label{xdifferenziabile}
  Let \ref{ipotesiuno} hold. Then, for each $(t,x)\in[0,T]\times H$, \eqref{eq:markf} has a unique solution $\{ X^{t,x}_\tau\}_{0\leq \tau\leq T}$. Moreover, for every $p>1$,

  \begin{description} 
    \item[(i)]   $X^{t,x}$ belongs to $\ys(H)$ and there exists a constant $C$ such that 
    \begin{equation}
    \label{stimadeimomentidix}
      \E\left[ \sup\nl _{\tau\in[0,T]}|X^{t,x}_\tau |^p \right]\leq C(1+|x|)^p,
    \end{equation}
    \item[(ii)]   The map
  $(t,x)\mapsto X^{t,x}$ belongs to $\calg^{0,1}\Big([0,T]\times
  H, \m S^p(H)\Big)$.

    \item[(iii)]  For every $h\in H$, the directional derivative process $\nabla_xX_{\t}^{t,x}h$,
 $\tau\in [0,T]$, solves the equation:
$$
\left\{
\begin{array}{lll}\dis
\nabla_xX_{\t}^{t,x}h&=&\dis
 e^{(\tau-t)A}h+\int_t^\tau e^{(\tau-r)A}
\nabla_xb(r,X_r^{t,x})\nabla_xX_r^{t,x}h\; dr
\\
&&+\dis
\int_t^\tau \nabla_x( e^{(\tau-r)A}
\sigma(r,X_r^{t,x}))\nabla_xX_r^{t,x}h\; dW_r,
\quad \tau\in [t,T],
\\\dis
\nabla_xX_{\t}^{t,x}h&=&h,\quad \tau\in [0,t).
\end{array}\right.
$$

\item[(iii)]  Finally
$\left\|\nabla_xX_{\t}^{t,x}h \right\|_{{\m S}^p}\leq c\,|h|$ for some constant $c$.
   \end{description}

\end{prop}

We assume that $F: [0,T]\times H \times\rset\times L_2(\Xi,\rset)\fl \rset$ and $\Phi:H\fl\rset$ are measurable functions such that

\begin{hyp}
\label{hdiff}
There exists $C\geq 0$ and $\alpha\in(0,1)$ such that
\begin{itemize}
\item $|F(t,x,y,z)| \leq C\left( 1+|y| + |z|^2 \right)$ and $\Phi$ is bounded ;
\item $F(s,\cdot,\cdot,\cdot)$ is $\mathcal{G}^{1,1,1}(H\times \R\times L_2(\Xi, \R);\R)$ and $\Phi$ is $\m G^1(H;\R)$ ;
\item  $|\nabla_x \Phi(x)| \leq C\left(1+|x|^n\right)$ ;
\item $\left| \nabla_x F(s,x,y,z) \right| \leq C\left( 1 + |x|^n + |z|^2 \right)$ ;
\item $\left| \nabla_z F(s,x,y,z) \right| \leq C\left( 1  + |z| \right)$ ;
\item $\left| \nabla_y F(s,x,y,z) \right| \leq C\left( 1  + |z| \right)^{2\alpha}$ ;
\end{itemize}
\end{hyp}

We know from results of \cite{Kob00, LsM98} (these results can be easily generalised to the case of a cylindrical Wiener process) that under $\ref{hdiff}$ the BSDE (\ref{eq:markb})
has a unique bounded solution and that there  exists a constant $C$ such that, for each $(t,x)$,
\begin{equation}
\label{equbmo}
\left\| \sup\nl_{u\in[0,T]} \left| Y_u^{t,x} \right| \right\|_\infty + \left\| \int_0^\cdot Z_s^{t,x}\cdot dW_s\right\|_{BMO_2} \leq C.
\end{equation} 
For the existence and the bound for the process $Y$ we refer to \cite[Corollary 1]{LsM98}, uniqueness follows from \cite[Theorem 2.6]{Kob00} and finally the estimate for the BMO-norm of $Z$ comes from a direct computation starting from It\^o's formula applied to $\varphi(x)= \left(e^{2Cx}-2Cx-1\right)/(2C^2)$.
In particular, for each $p\geq 1$,
\begin{equation}
\label{eqz2}
\left\| \Big(\int_0^T \left| Z_s^{t,x} \right|^2 ds\Big)^{1/2} \right\|_p \leq C_p.
\end{equation}

\begin{prop}\label{ydiff}
Let the assumption \ref{hdiff} hold. 

The map $(t,x)\longmapsto \left( Y_{\cdot}^{t,x}, Z_{\cdot}^{t,x} \right)$ belongs to $\m G^{0,1}\left([0,T]\times H ; \ys\times\zs \right)$ for each $p>1$. Moreover, for every $x\in H$ and $h\in H$,
  the directional derivative process 
$\left\{ \nabla_xY_u^{t,x}h, \nabla_xZ_u^{t,x}h \right\}_{u\in [0,T]}$ solves the
BSDE: for $\tau\in[0,T]$,
\begin{equation}
\label{eqdp}
\begin{split}
\nabla_x Y^{t,x}_u  h = & \nabla_x \Phi \left(X^{t,x}_T \right) \nabla_x X^{t,x}_T h +  \int_u^T \nabla_x F\left(s,X^{t,x}_s,Y^{t,x}_s,Z^{t,x}_s\right) \nabla_x X^{t,x}_s h \, ds \\
 &\hspace*{-10mm} + \int_u^T \left( \nabla_y F\left(s,X^{t,x}_s,Y^{t,x}_s,Z^{t,x}_s\right) \nabla_x Y^{t,x}_s h + \nabla_z F\left(s,X^{t,x}_s,Y^{t,x}_s,Z^{t,x}_s\right) \nabla_x Z^{t,x}_s h \right) ds \\
 & - \int_u^T \nabla_x Z^{t,x}_s h\, dW_s
\end{split}
\end{equation}
and there exists $C_p$ such that
$$
\left\| \nabla_x Y^{t,x}h \right\|_{\ys} + \left\| \nabla_x Z^{t,x} h \right\|_{\zs} \leq C_p(1+|x|)^n|h|.
$$

\end{prop}

\begin{proof} 
    The continuity of the map $(t,x)\longmapsto \left( Y_{\cdot}^{t,x}, Z_{\cdot}^{t,x} \right)$ follows from a mere extension of Kobylanski's stability result \cite[Theorem 2.8]{Kob00}.

For the differentiability, let us remark that, in view of \ref{hdiff} and \eqref{eqz2}, for all $p>1$,
$$\left\| \left| \nabla_x \Phi\left(X^{t,x}_u \right) \nabla_x X^{t,x}_T h \right| + \int_0^T \left| \nabla_x F\left(s,X^{t,x}_s,Y^{t,x}_s,Z^{t,x}_s\right) \nabla_x X^{t,x}_s h \right| ds \right\|_p \leq C_p (1+|x|)^n|h|.
$$
It follows from Theorem~\ref{resesun}, that the BSDE~\eqref{eqdp} has a unique solution which belongs to $\ys \times \zs$ for all $p\geq 1$. And moreover, for $p>1$, it follows from Corollary~\ref{resutil} and \eqref{eqz2}, that
$$
\left\| \nabla_x Y^{t,x} h \right\|_{\ys} + \left\| \nabla_x Z^{t,x} h \right\|_{\zs} \leq C(1+|x|)^n|h|.
$$

Let us fix $(t,x)\in[0,T]\times H$. We remove the parameters $t$ and $x$ for notational simplicity. For $\ep>0$, we set $X^\ep = X^{t,x+\ep h}$, where $h$ is some vector in $H$, and we consider $(Y^\ep,Z^\ep)$ the solution in $\m S^p \times M^p$ to the BSDE
$$
Y^\ep_t = \Phi(X_T^{t,\ep})+ \int_t^T F(s,X^\ep_s,Y^\ep_s,Z^\ep_s)\,ds - \int_t^T Z^\ep_s\,dW_s.
$$
When $\ep\to 0$, $\left(X^\ep,Y^\ep,Z^\ep\right)\fl (X,Y,Z)$ in $\ys\times\ys\times\zs$ for all $p>1$. We also denote  $(G,N)$ the solution to the BSDE~\eqref{eqdp} and it remains to prove that the directional derivative of the map $(t,x)\longmapsto \left( Y_{\cdot}^{t,x}, Z_{\cdot}^{t,x} \right)$ in the direction $h \in H$ is given by $(G,N)$.

Let us consider $U^\ep = \ep^{-1}\left(Y^\ep-Y\right) - G$, $V^\ep = \ep^{-1}\left(Z^\ep-Z\right) - N$.
 We have,
\begin{eqnarray*}
U^\ep_t  & = & \frac{1}{\ep}\left(\Phi(X^\ep_T) - \Phi(X_T)\right) - \nabla_x \Phi(X_T) \nabla_x X_T h +
\\
& & +\frac{1}{\ep} \int_t^T \left(F(s,X^\ep_s,Y^\ep_s,Z^\ep_s) - F(s,X_s,Y_s,Z_s)\right) ds  - \int_t^T V^\ep_s\,dW_s  \\
& & - \int_t^T \nabla_x F(s,X_s,Y_s,Z_s) \nabla_x X_s h \, ds - \int_t^T \nabla_y F(s,X_s,Y_s,Z_s) G_s \, ds \\
& & - \int_t^T \nabla_z F(s,X_s,Y_s,Z_s) N_s \, ds.
\end{eqnarray*}
Using the fact that $\psi(s,\cdot,\cdot,\cdot)$ belongs to $\m G^{1,1,1}$, we can write
$$
\frac{1}{\ep}\left(F(s,X^\ep_s,Y^\ep_s,Z^\ep_s) - F(s,X_s,Y_s,Z_s)\right) = \frac{1}{\ep}\left(F(s,X^\ep_s,Y_s,Z_s) - F(s,X_s,Y_s,Z_s)\right)+$$
$$ \quad \quad \quad + A^\ep_{s} \, \frac{Y^\ep_{s}-Y_s}{\ep} + B^\ep_s\, \frac{Z^\ep_s - Z_s}{\ep}
$$
where $A^{\ep}_s \in L(\R, \R)$ and $B^\ep_s\in L\left( L_2(\Xi,\R), \R \right)$ are defined by
$$
\forall y\in \R,\qquad A^\ep_s y = \int_0^1 \nabla_y F\left(s,X^\ep_s,Y_s +\alpha (Y^\ep_s-Y_s), Z_s\right) y \,d\alpha,
$$
$$
\forall z\in L_2(\Xi,\R),\qquad B^\ep_s z = \int_0^1 \nabla_z F\left(s,X^\ep_s,Y^\ep_s,Z_s +\alpha (Z^\ep_s-Z_s)\right) z \,d\alpha.
$$

Then $(U^\ep,V^\ep)$ solves the following BSDE
$$
U^\ep_t = \zeta^\ep +\int_t^T \left( A^\ep_s U^\ep_s + B^\ep_s V^\ep_s\right) ds + \int_t^T \left( P^\ep(s) + Q^\ep(s) + R^\ep(s) \right) ds - \int_t^T V^\ep_s\,dW_s
$$
where we have set
$$
P^\ep(s) = \left( A^\ep_s - \nabla_y F(s,X_s,Y_s,Z_s) \right) G_s, \qquad Q^\ep(s) = \left( B^\ep_s - \nabla_z F(s,X_s,Y_s,Z_s) \right) N_s,
$$
$$
R^\ep(s) = \ep^{-1}\left( F(s,X^\ep_s,Y_s,Z_s) - F(s,X_s,Y_s,Z_s)\right) - \nabla_x F(s,X_s,Y_s,Z_s) \nabla_x X_s h,
$$
$$
\zeta^\ep = \ep^{-1}\left( \Phi(X^\ep_T) - \Phi(X_T)\right) - \nabla_x \Phi(X_T) \nabla_x X_T h.
$$

It follows from \ref{hdiff} that
$$
A^\ep_s \leq C\left( 1 + |Z_s| + \left| Z^\ep_s\right| \right)^{2\alpha}, \qquad \left| B^\ep_s\right| \leq C\left( 1 + |Z_s| + \left| Z^\ep_s\right| \right),
$$
and
$$
\left| P^\ep(s) \right| \leq C \left( 1 + |Z_s| + \left| Z^\ep_s\right| \right)^{2\alpha} |G_s|, \qquad \left| Q^\ep(s) \right| \leq C \left( 1 + |Z_s| + \left| Z^\ep_s\right| \right) |H_s|
$$
For $p$ large enough, we have from Corollary~\ref{resutil} taking into account \eqref{equbmo} and \eqref{eqz2},
$$
\left\| U^\ep\right\|_{\ys} + \left\| V^\ep \right\|_{\zs}  \leq C\, \left\| \left|\zeta^\ep \right| +  \int_0^T \left( \left| P^\ep(s) \right| + \left| Q^\ep(s) \right| + \left| R^\ep(s) \right|\right)ds \right\|_{p+1}.
$$
The right hand side of the previous inequality tends to 0 as $\ep\to 0$  in view of the regularity and the growth of $F$ and $\Phi$ (see \ref{hdiff}).

The proof that the maps $x \mapsto (\nabla_x Y^{t,x} h, \nabla_x Z^{t,x}h)$ and $h \mapsto (\nabla_x Y^{t,x}h, \nabla_x Z^{t,x}h)$ are continuous (for every $h$ and $x$ respectively) comes once again of Corollary~\ref{resutil}.
\end{proof}

\begin{rem}
Since $\sup_{t,x} \left\| \sup_u \left| Y(u,t,x)\right|\right\|_\infty <\infty$, one can change $C$ by $C(|y|)$ in the assumptions on the gradient on $F$ in \ref{hdiff}.
\end{rem}

\section{Application to nonlinear PDEs}
\label{sec-nlpdes}
In this section we are interested in finding a probabilistic representation in our framework for the solution to 
\begin{equation}
\label{kolmogorovnonlineare}
  \left\{\begin{array}{l}\dis
\partial_t u(t,x) +\call_t [u(t,\cdot)](x) +
F (t, x,u(t,x),\sigma(t,x)^*\nabla_xu(t,x)) = 0,\quad t\in [0,T],\,
x\in H,\\
\dis u(T,x)=\Phi(x),
\end{array}\right.
\end{equation} 
where $\call_t$ is the operator:
$$
\call_t[\phi](x)=\frac{1}{2}{\rm Trace }\left(
\sigma(t,x)\sigma(t,x)^*\nabla^2\phi(x)\right) + \< Ax+b(t,x),\nabla\phi(x)\>,
$$
where $\nabla\phi$ and $\nabla^2\phi$
are the first and the second G\^ateaux derivatives of
$\phi$ (identified with elements of $H$ and $L(H)$ respectively).
This definition is formal, since the domain of $\call_t$ is
not specified.  

We will refer to this equation as the nonlinear Kolmogorov
equation. In this equation, $F: [0,T]\times H\times \R\times \Xi\to \R$ is a given
function verifying \ref{hdiff} and $\nabla_xu(t,x)$ is the
G\^ateaux derivative of
$u(t,x)$ with respect to $x$: it is identified with an element
of $H$, so that $\sigma(t,x)^*\nabla_xu(t,x)\in \Xi$.

\smallskip

Under the assumption \ref{ipotesiuno}, we can define a transition semigroup $P_{t,\tau}$ with the help of  $X^{t,x}$ solution to \eqref{eq:markf} by the formula
$$
P_{t,\tau}[\phi](x)=\E\left[ \phi(X_{\t}^{t,x}) \right],\qquad x\in H.
$$
The estimate \eqref{stimadeimomentidix} shows that
$P_{t,\tau}$ is well defined as a linear operator from $
\calb_p(H)$, the set of measurable functions from $H$ to $\rset$ with polynomial growth, into itself; the semigroup property
$P_{t,s}P_{s,\tau}=P_{t,\tau}$, $t\leq s\leq \tau$,
is well known.

When $\phi$ is sufficiently regular, the function
$v(t,x)=P_{t,T}[\phi](x)$, is a classical solution of the backward
Kolmogorov equation \eqref{kolmogorovnonlineare} with $F\equiv 0$;
 we refer to \cite{dPZ92} and \cite{Zab99} for a detailed
exposition.
When $\phi$ is not regular,
the function $v$ defined by the formula
$v(t,x)=P_{t,T}[\phi](x)$ can be considered as a
generalized solution of this equation.

For the nonlinear case, we consider  the variation of
constants formula for
 (\ref{kolmogorovnonlineare}):
\begin{equation}
\label{solmild}
  u(t,x) =\int_t^TP_{t,\tau}[F (\tau, \cdot,u(\tau,\cdot),\sigma(\tau,\cdot)^*
\nabla_xu(\tau,\cdot))
](x)\; d\tau+ P_{t,T}[\Phi](x),\quad t\in [0,T],\,
x\in H,
\end{equation}
and we  notice that this formula  is meaningful,
provided $F(t,\cdot,\cdot,\cdot)$,
$u(t,\cdot)$ and $\nabla_xu(t,\cdot)$ have polynomial
growth.
We use this formula as a definition for the solution of
(\ref{kolmogorovnonlineare}):
\begin{df}\label{defdisoluzionemild}
We say that a function
$u:[0,T]\times H\to\R$ is a mild solution of the nonlinear
Kolmogorov equation
(\ref{kolmogorovnonlineare}) if the following conditions hold:
\begin{description}
  \item[(i)]
$u\in\calg^{0,1}([0,T]\times H,\R)$;
  \item[(ii)]  there exists $C>0$ and $d\in \N$ such that
  $|\nabla_xu(t,x)h|\leq C|h|(1+|x|^{d})$ for all $t\in [0,T]$, $x\in H$,
$h\in H$;
 \item[(iii)] equality (\ref{solmild}) holds.
\end{description}
\end{df}

\begin{rem}
We obtain an equivalent formulation
of (\ref{kolmogorovnonlineare})
and  (\ref{solmild})
by considering
the G\^ateaux derivative $\nabla_xu(t,x)$ as an element of
$\Xi^*= L(\Xi,\R)=L_2(\Xi,\R)$. In this case, we  take a function
$F: [0,T]\times H\times \R\times L_2(\Xi,\R)\to \R$ and we
write the equation in the form
$$
\partial_t u(t,x)+\call_t [u(t,\cdot)](x) +
F (t, x,u(t,x),\nabla_xu(t,x)\sigma(t,x)) =0.
$$
The two forms are equivalent provided we identify
$\Xi^*= L_2(\Xi,\R)$ with $\Xi$ by the Riesz isometry.
\end{rem}

We are now ready to state the main result of this section.
\begin{thm}\label{main}
Let the assumptions \ref{ipotesiuno} and \ref{hdiff} hold.

The nonlinear Kolmogorov equation (\ref{kolmogorovnonlineare}) has a unique mild solution $u$ given by the formula
$$
u(t,x) =Y^{t,x}_t,\quad (t,x)\in[0,T]\times H
$$
where $\left(Y^{t,x},Z^{t,x}\right)$ is the solution to the BSDE~\eqref{eq:markb} and $X^{t,x}$ the solution to \eqref{eq:markf}. Moreover, we have, $\P$--a.s.
$$
Y_s^{t,x} = u(s,X_s^{t,x}),\qquad Z_s^{t,x}\sigma(s, X_s^{t,x})^*\nabla_xu(s,t, X_s^{t,x}).
$$
\end{thm}

\begin{proof}
Let us first recall a result of \cite[Lemma 6.3]{FT02}.
Let $\{ e_i\}$ be a basis of  $\Xi$ and let us consider
the standard real Wiener process $W^i_\tau = \int_0^\tau \< e_i,dW_\sigma\>$, $\tau\geq 0$.

If $v\in \m G^{0,1}([0,T]\times H, \rset)$, for every $i$, the  quadratic variation
  of $ v(s,X_s^{t,x})$ and $W^i_s$ is given by
\begin{equation}
\label{quadrvardiu}
\left[v(\cdot,X_{\cdot}^{t,x}),W^i\right]_s = \int_t^s \nabla_x v(\tau, X_{\t}^{t,x})G(\tau, X_{\t}^{t,x})e_i
\; d\tau,\quad s\in [t,T].
\end{equation}

{\em (a) Existence.}
Let us recall that for $s\in[t,T]$, $Y_s^{t,x}$ is measurable with respect to $\m F_{[t,s]}$ and $\m F_s$; it follows that $Y_t^{t,x}$ is deterministic (see also \cite{ElK97}). Moreover, as a byproduct of Proposition~\ref{ydiff}, the function $u$ defined by the formula $u(t,x)=Y^{t,x}_t$ has the
regularity properties stated in Definition~\ref{defdisoluzionemild}. It
remains to verify that
equality (\ref{solmild}) holds true for $u$.

To this purpose we first fix $t\in [0,T]$ and $x\in H$. Since $(Y_{\cdot}^{t,x},Z_{\cdot}^{t,x})$ solves the BSDE \eqref{eq:markb}, we have, for $s\in[t,T]$,
$$
 Y_s^{t,x}+\int_s^TZ_{\t}^{t,x}\, dW_\tau= \Phi(X_T^{t,x}) +
  \int_s^T F \Big(\tau ,X_{\t}^{t,x},
   Y_{\t}^{t,x},Z_{\t}^{t,x}\Big)\, d\tau,
$$
and, taking expectation for $s=t$ we obtain, coming back to the definition of $u$ and $P_{t,T}$,
\begin{equation}
\label{presque}
 u(t,x)= P_{t,T}[\Phi](x) +
  \E\left[\int_t^T F \Big(\tau,X_{\t}^{t,x}, Y_{\t}^{t,x},Z_{\t}^{t,x}
  \Big)\; d\tau \right].
\end{equation}
Moreover, we have, for each $i$,
$$
\left[ Y_{\cdot}^{t,x},W^i\right]_s = \int_t^s \<Z_\tau, e_i \> \,d\tau, \quad s\in[t,T].
$$

Now let us observe that the processes $Y$ and $Z$ satisfy the Markov property: for $t \leq s\leq T$, $\P$-a.s.
$$Y_{\t}^{s,X_s^{t,x}}=Y_{\t}^{t,x} \quad \textit{for } \t \, \in [s,T]$$
$$ Z_{\t}^{s,X_s^{t,x}}=Z_{\t}^{t,x} \quad \textit{for a.e. } \t \, \in [s,T].$$

In fact the solution of the backward equation is uniquely determined on an interval $[s,T]$ by the values of the process $X$ on the same interval. The process $X$ is the unique solution of the forward equation (\ref{eq:markf}) and satisfies the Markov property.

As consequence we have, $\P$--a.s.,
$$
u(\t,X_{\t}^{t,x})= Y_{\t}^{t,x}, \quad \t \in [t,T].
$$
It follows from \eqref{quadrvardiu} that, for each $i$,
$$
\left[ Y_{\cdot}^{t,x},W^i\right]_s = \int_t^s \nabla_x u(\tau, X_{\t}^{t,x})\sigma(\tau, X_{\t}^{t,x})e_i\, d\tau,\quad s\in [t,T].
$$
Therefore, for a.a. $\tau\in [t,T]$, we have $\P$-a.s.
$$
\nabla_xu(\tau, X_{\t}^{t,x})\sigma(\tau, X_{\t}^{t,x})e_i = \<Z_{\t}^{t,x},e_i\>,
$$
for every $i$. Identifying $\nabla_x u(t,x)$ with an element of
$\Xi$, we conclude that for a.a. $\tau\in [t,T]$,
$$
\sigma(\tau, X_{\t}^{t,x})^*\nabla_xu(\tau,t, X_{\t}^{t,x}) = Z_{\t}^{t,x}.
$$

Thus, $ F \left(\tau,X_{\t}^{t,x}, Y_{\t}^{t,x},Z_{\t}^{t,x}
  \right) $ can be rewritten as
$$
F \left(\tau, X_{\t}^{t,x},u(\tau,X_{\t}^{t,x}),\sigma(\tau,X_{\t}^{t,x})^*
\nabla_xu(\tau,X_{\t}^{t,x})\right)
$$
and \eqref{presque} leads to
$$
u(t,x)= P_{t,T}[\phi](x) + \int_t^TP_{t,\tau}[
F (\tau, \cdot,u(\tau,\cdot),\sigma(\tau,\cdot)^*
\nabla_xu(\tau,\cdot))
](x)\, d\tau
$$
which is \eqref{solmild}.

\medskip

{\em (b) Uniqueness.}
Let $u$ be a mild solution. We look for a convenient
expression for the process
$u(s,X_s^{t,x})$, $s\in [t,T]$.
By (\ref{solmild}) and the definition of $P_{t,\tau}$,
for every $s\in [t,T]$ and $x\in H$,
\begin{eqnarray*}
  u(s,x)  & = & \E \left[\Phi(X_T^{s,x}) \right] \\
 & & +\E\left[\int_s^T F\big(\tau, X_{\t}^{t,x},u(\tau,X_{\t}^{t,x}),
\sigma(\tau,X_{\t}^{t,x})^*
\nabla_xu(\tau,X_{\t}^{t,x})\big)
d\tau\right].
\end{eqnarray*}

Since $X_{\t}^{t,x}$ is independent of $\calf_s$, we can replace the
expectation by the conditional expectation given $\calf_s$:
\begin{eqnarray*}
  u(s,x)  & = & \E^{\calf_s} \left[\Phi(X_T^{s,x}) \right] \\
 & & +\E^{\calf_s}\left[\int_s^T F\big(\tau, X_{\t}^{t,x},u(\tau,X_{\t}^{t,x}),
\sigma(\tau,X_{\t}^{t,x})^*
\nabla_xu(\tau,X_{\t}^{t,x})\big)
d\tau\right].
\end{eqnarray*}

Taking into account the Markov property of $X$, $\P$--a.s.
$$
X_{\tau}^{s,X_s^{t,x}}=X_{\t}^{t,x},
\qquad \tau\in [s,T],
$$
we have
\begin{eqnarray*}
 \lefteqn{u(s,X_s^{t,x}) =\E^{\calf_s} \left[ \Phi(X_T^{t,x}) \right]  } \\
   & & \qquad
+\E^{\calf_s}\left[\int_s^T
F\big(\tau, X_{\t}^{t,x},u(\tau,X_{\t}^{t,x},
\sigma(\tau,X_{\t}^{t,x})^*
\nabla_xu(\tau,X_{\t}^{t,x}) \big)
d\tau\right].
\end{eqnarray*}
If we set
$$
\xi=\Phi(X_T^{t,x})+\int_t^T
F\big(\tau, X_{\t}^{t,x},u(\tau,X_{\t}^{t,x}),
\sigma(\tau,X_{\t}^{t,x})^*
\nabla_xu(\tau,X_{\t}^{t,x})\big) d\tau
$$
the previous equality leads to
\begin{eqnarray*}
\lefteqn{u(s,X_s^{t,x})} \\
&&  =\E^{\calf_s} \,[\xi] -\int_t^s
F \big(\tau, X_{\t}^{t,x},u(\tau,X_{\t}^{t,x}),
\sigma(\tau,X_{\t}^{t,x})^*
\nabla_xu(\tau,X_{\t}^{t,x})\big)
 d\tau.
\end{eqnarray*}

Let us observe that $\E^{\calf_t}[\xi] = u(t,x)$.
Since $\xi\in L^2(\Omega;\R)$ is $\calf_{[t,T]}$--measurable,
by the representation theorem, there exists
$\widetilde{Z}\in L^2_\calp(\Omega \times [t,T]; L_2(\Xi,\R))$ such
that
$$
\E^{\calf_s} [\xi] =u(t,x) + \int_t^s\widetilde{Z}_\tau\, dW_\tau,\quad s\in[t,T].
$$
We conclude that
 the process
$u(s,X_s^{t,x})$, $s\in [t,T]$ is a (real) continuous
semimartingale with canonical decomposition
\begin{eqnarray}
\label{scomposizione}
u(s,X_s^{t,x}) & = & u(t,x) + \int_t^s\widetilde{Z}_\tau\; dW_\tau\\
\nonumber
 &&-\int_t^s
F \Big(\tau, X_{\t}^{t,x},u(\tau,X_{\t}^{t,x}),
\sigma(\tau,X_{\t}^{t,x})^*
\nabla_xu(\tau,X_{\t}^{t,x})\Big)
\, d\tau.
\end{eqnarray}

Using \eqref{quadrvardiu} and arguing as in the proof of existence, we deduce that
for a.a. $\tau\in [t,T]$, $\P$-a.s.
$$
\sigma(\tau, X_{\t}^{t,x})^*\nabla_xu(\tau, X_{\t}^{t,x})=\widetilde{Z}_\tau.
$$
Substituting into \eqref{scomposizione} we obtain
\begin{eqnarray*}
u(s,X_s^{t,x}) & = &  u(t,x) +
\int_t^s \sigma(\tau, X_{\t}^{t,x})^*\nabla_xu(\tau, X_{\t}^{t,x})\, dW_\tau \\
&& -\int_t^s
F \big(\tau, X_{\t}^{t,x},u(\tau,X_{\t}^{t,x}),
\sigma(\tau,X_{\t}^{t,x})^*
\nabla_xu(\tau,X_{\t}^{t,x})\big)
\, d\tau,
\end{eqnarray*}
for $s\in [t,T]$. Since $u(T,X_T^{t,x})=\Phi(X_T^{t,x})$, we deduce that
$$
\left\{ \big(u(s,X_s^{t,x}, \sigma(\tau,X_{\t}^{t,x})^*\nabla_x u(\tau,X_{\t}^{t,x}\big) \right\}_{s\in[t,T]}
$$
solves
the backward equation \eqref{eq:markb}. By uniqueness, we have $Y_s^{t,x}=u(s,X_s^{t,x})$, for each $s\in [t,T]$ and in particular, for $s=t$, $u(t,x)=Y_t^{t,x}$.

\end{proof}

\end{document}